\documentclass[a4paper,12pt]{amsart}
\usepackage{amsfonts}
\usepackage{amssymb}
\usepackage{ifthen}
\usepackage{graphicx}
\nonstopmode \numberwithin{equation}{section}
\usepackage[usenames]{color}
\newtheorem{thm}{Theorem}[section]
\newtheorem{lem}{Lemma}[section]
\newtheorem{cor}{Corollary}[section]

\newtheorem{cl}{Claim}[section]
\newtheorem{ca}{Case}[section]
\newtheorem{sca}{Subcase}[section]
\newtheorem{scl}[section]{Subclaim}
\newtheorem{conj}[equation]{Conjecture}

\theoremstyle{definition}
\newtheorem{defn}{Definition}[section]
\newtheorem{op}[equation]{Open Problem}
\newtheorem{ques}[equation]{Question}
\newtheorem{rem}{Remark}[section]
\newtheorem{exam}[equation]{Example}

\newcounter {own}
\def\theown {\thesection       .\arabic{own}}

\newenvironment{pf}[1][]{%
 \vskip 3mm
 \noindent
 \ifthenelse{\equal{#1}{}}%
  {{\slshape Proof. }}%
  {{\slshape #1.} }%
 }%
{\qed\bigskip}

\newcounter{alphabet}
\newcounter{tmp}
\newenvironment{Thm}[1][]{\refstepcounter{alphabet}%
\bigskip%
\noindent%
{\bf Theorem \Alph{alphabet}}%
\ifthenelse{\equal{#1}{}}{}{ (#1)}%
{\bf .} \itshape}{\vskip 8pt}

\makeatletter
\newcommand{\Ref}[1]{\@ifundefined{r@#1}{}{\setcounter{tmp}{\ref{#1}}\Alph{tmp}}}
\makeatother

\newenvironment{Lem}[1][]{\refstepcounter{alphabet}%
\bigskip%
\noindent%
{\bf Lemma \Alph{alphabet}}%
{\bf .} \itshape}{\vskip 8pt}

\newcommand{\IR}{{\mathbb R}}

\newcommand{\diam}{{\operatorname{diam}}}

\newcommand{\dist}{{\operatorname{dist}}}

\def\be{\begin{equation}}
\def\ee{\end{equation}}

\newcommand{\ben}{\begin{enumerate}}
\newcommand{\een}{\end{enumerate}}

\newcommand{\blem}{\begin{lem}}
\newcommand{\elem}{\end{lem}}
\newcommand{\bthm}{\begin{thm}}
\newcommand{\ethm}{\end{thm}}
\newcommand{\bcor}{\begin{cor}}
\newcommand{\ecor}{\end{cor}}
\newcommand{\beg}{\begin{exam}}
\newcommand{\eeg}{\end{exam}}
\newcommand{\begs}{\begin{examples}}
\newcommand{\eegs}{\end{examples}}
\newcommand{\bdefe}{\begin{defn}}
\newcommand{\edefe}{\end{defn}}
\newcommand{\bprob}{\begin{prob}}
\newcommand{\eprob}{\end{prob}}
\newcommand{\bques}{\begin{ques}}
\newcommand{\eques}{\end{ques}}
\newcommand{\bei}{\begin{itemize}}
\newcommand{\eei}{\end{itemize}}
\newcommand{\bcon}{\begin{conj}}
\newcommand{\econ}{\end{conj}}
\newcommand{\bop}{\begin{op}}
\newcommand{\eop}{\end{op}}

\newcommand{\bca}{\begin{ca}}
\newcommand{\eca}{\end{ca}}
\newcommand{\bsca}{\begin{sca}}
\newcommand{\esca}{\end{sca}}

\newcommand{\bcl}{\begin{cl}}
\newcommand{\ecl}{\end{cl}}

\newcommand{\bscl}{\begin{scl}}
\newcommand{\escl}{\end{scl}}

\newcommand{\bcons}{\begin{conjs}}
\newcommand{\econs}{\end{conjs}}
\newcommand{\bprop}{\begin{propo}}
\newcommand{\eprop}{\end{propo}}
\newcommand{\br}{\begin{rem}}
\newcommand{\er}{\end{rem}}
\newcommand{\brs}{\begin{rems}}
\newcommand{\ers}{\end{rems}}
\newcommand{\bo}{\begin{obser}}
\newcommand{\eo}{\end{obser}}
\newcommand{\bos}{\begin{obsers}}
\newcommand{\eos}{\end{obsers}}
\newcommand{\bpf}{\begin{pf}}
\newcommand{\epf}{\end{pf}}
\newcommand{\ba}{\begin{array}}
\newcommand{\ea}{\end{array}}
\newcommand{\beq}{\begin{eqnarray}}
\newcommand{\beqq}{\begin{eqnarray*}}
\newcommand{\eeq}{\end{eqnarray}}
\newcommand{\eeqq}{\end{eqnarray*}}

\newcounter{minutes}
\divide\time by 60
\newcounter{hours}
\multiply\time by 60 \addtocounter{minutes}{-\time}

\begin{document}

\bibliographystyle{amsplain}
\title{Gromov hyperbolicity, John spaces and quasihyperbolic geodesics}

\author{Qingshan Zhou}
\address{Qingshan Zhou, School of Mathematics and Big Data, Foshan University,  Foshan, Guangdong 528000, People's Republic
of China} \email{q476308142@qq.com}

\author{Yaxiang  Li${}^{\mathbf{*}}$}
\address{Yaxiang Li,  Department of Mathematics, Hunan First Normal University, Changsha,
Hunan 410205, People's Republic
of China} \email{yaxiangli@163.com}

\author{Antti Rasila}
\address{Antti Rasila, College of Science, Guangdong Technion -- Israel Institute of Technology, Shantou, Guangdong 515063, People's Republic
of China} \email{antti.rasila@gtiit.edu.cn; antti.rasila@iki.fi}

\def\thefootnote{}
\footnotetext{ \texttt{\tiny File:~\jobname .tex,
          printed: \number\year-\number\month-\number\day,
          \thehours.\ifnum\theminutes<10{0}\fi\theminutes}
} \makeatletter\def\thefootnote{\@arabic\c@footnote}\makeatother

\date{}
 \keywords{ Quasihyperbolic metric, Gromov hyperbolic spaces, John spaces, quasihyperbolic geodesic.\\
${}^{\mathbf{*}}$ Corresponding author}

\begin{abstract} We show that every quasihyperbolic geodesic in a John space admitting a roughly starlike Gromov hyperbolic quasihyperbolization is a cone arc. This result provides a new approach to the elementary metric geometry question, formulated in \cite[Question 2]{Hei89}, which has been studied by Gehring, Hag, Martio and Heinonen. As an application, we obtain a simple geometric condition connecting uniformity of the space with the existence of Gromov hyperbolic quasihyperbolization.
\end{abstract}

\thanks{The research was partly supported by NNSF of
China (Nos. 11601529,  11671127,  11571216).}

\maketitle{} \pagestyle{myheadings} \markboth{}{}

\section{Introduction}
The unit disk or Poincar\'e disk $\mathbb{D}$ serves as a canonical model in studying of conformal mappings and hyperbolic geometry in complex analysis. It  is  noncomplete metric space with the metric inherited from the two dimensional Euclidean space $\mathbb{R}^2$. On the other hand, the unit disk equipped with the Poincar\'e metric is  complete Riemannian $2$-manifold with constant negative curvature. This observation can be used in investigating the hyperbolic metric on planar domains and conformal mappings between them. A generalization of this idea to higher dimensional spaces, involving quasihyperbolic metrics and Gromov hyperbolicity, was studied by Bonk, Heinonen and Koskela in \cite{BHK}.

Well-known geometric properties of a hyperbolic geodesic $[x,y]\in \mathbb{D}$ with respect to the Euclidean metric are:
\begin{itemize}
  \item $\ell([x,y])\leq C|x-y|$,
  \item $\min\{\ell([x,z]),\ell([z,y])\}\leq C\dist(z,\partial \mathbb{D})$
\end{itemize}
for all $z\in [x,y]$, where $C$ is a universal constant. The first of the above conditions says that hyperbolic geodesic essentially minimizes the length of all curves connecting the endpoints, namely, the Gehring-Haymann condition. The second one is called the cone condition or the double twisted condition.

Martio and Sarvas studied in \cite{MS78} global injectivity properties of  locally injective mappings. They considered a class of domains of $\mathbb{R}^n$, named by {\it uniform domains},  which means every pair of points can be connected by a curve satisfies the above two conditions for some constant $C\geq 1$. In \cite{GO}, Gehring and Osgood investigated the geometric properties of {\it quasihyperbolic metric}, which was introduced by Gehring and Palka \cite{GP76}, and proved that every quasihyperbolic geodesic in a Euclidean uniform domain also satisfies the above two conditions.

It should be noted that the class of domains on $\mathbb{R}^n$, which only satisfies the second condition known as {\it John domains} is large and of independent interest. For instance, the slit disk on $\mathbb{R}^2$ is an example of such domain. This class was first considered by John \cite{Jo61} in the context of elasticity theory. Many characterizations of uniform and John domains can be found in the literature and the importance of these classes of domains in function theory is well established, see for example \cite{GGKN17, LVZ17}.


From a geometric point of view, it is natural question, whether each quasihyperbolic geodesic of a John domain is a cone arc. This question was pointed out already in 1989 by Gehring, Hag and Martio \cite{GHM}:

\bques\label{q-1}
Suppose $D\subset \mathbb{R}^n$ is a $c$-John domain and that $\gamma$ is a quasihyperbolic geodesic in $D$. Is $\gamma$ a $b$-cone arc for some $b=b(c)$?
\eques

They proved in \cite[Theorem $4.1$]{GHM} that quasihyperbolic geodesic in a plane simply connected John domain is a cone arc. They also constructed several examples to show that a similar result does not hold in higher dimensions. Furthermore, Heinonen has posed the following closely related problem concerning John disks:

\bques\label{q-2}$($\cite[Question 2]{Hei89}$)$
Suppose $D\subset \mathbb{R}^n$ is a $c$-John domain which is quasiconformally equivalent to the unit ball $\mathbb{B}$ and that $\gamma$ is a quasihyperbolic geodesic in $D$. Is $\gamma$ a $b$-cone arc for some constant $b$?
\eques

With the help of the conformal modulus of path families and Ahlfors $n$-regularity of $n$-dimensional Hausdorff measure of $\mathbb{R}^n$, Bonk, Heinonen and Koskela \cite[Theorem $7.12$]{BHK} gave an affirmative answer to Question \ref{q-2} for bounded domains with the constant dependence of the space dimension $n$. Recently, Guo \cite[Remark 3.10]{Guo15} provided a geometric method to deal with this question. His method was based on the result that a noncomplete metric space with a roughly starlike Gromov hyperbolic quasihyperbolization satisfies the Gehring-Hayman condition and the ball separation condition. These properties were established by Koskela, Lammi and Manojlovi\'{c} in \cite[Theorem $1.2$]{KLM14}. The constant $b$ in their results depends on the dimension $n$ as well.  The second author of this paper considered a related question for quasihyperbolic quasigeodesics  in the setting of Banach spaces \cite{Li}. Note that quasihyperbolic geodesics may not exist in infinite-dimensional spaces, even with assumption of convexity \cite{RT2}.

The concept of uniformity in a metric space setting was first introduced by Bonk, Heinonen and Koskela \cite{BHK}, where they connected the uniformity to the negative curvature of the space that is understood in the sense of Gromov. Moreover, they generalized the result of Gehring and Osgood and showed that every quasihyperbolic geodesic in a $c$-uniform space must be a $C$-uniform arc with $C=C(c)$, see \cite[Theorem 2.10]{BHK}. They also proved that $c$-uniform space is a Gromov $\delta$-hyperbolic with respect to its quasihyperbolic metric for some constant $\delta=\delta(c)$, see \cite[Theorem 3.6]{BHK}.

In view of the above results, it is natural to consider the following more general question:

\bques\label{q-3}
Let $D$ be a locally compact, rectifiably connected noncomplete metric space. If $D$ is an $a$-John space and $(D,k)$ is $\delta$-hyperbolic, is every quasihyperbolic geodesic $\gamma$ a $b$-cone arc with $b$ depending only on $a$ and $\delta$?
\eques

In this paper, we study these questions. Our main result is the following:

\begin{thm}\label{thm-1} Let $D$ be a locally compact, rectifiably connected noncomplete metric space. If $D$ is $a$-John and $(D,k)$ is $K$-roughly starlike and $\delta$-hyperbolic, then every quasihyperbolic geodesic in $D$ is a $b$-cone arc where $b$ depends only on $a, \delta$ and $K$.
\end{thm}
Every proper domain $D$ in $\mathbb{R}^n$ is a locally compact, rectifiably connected noncomplete metric space. Following terminology of \cite{BB03}, we call a locally compact, rectifiably connected noncomplete metric space $(D,d)$ {\it minimally nice}.
For a minimally nice space $(D,d)$, we say that $D$ has a {\it Gromov hyperbolic quasihyperbolization}, if $(D,k)$ is $\delta$-hyperbolic for some constant $\delta\geq 0$, where $k$ is the quasihyperbolic metric (for definition see Subsection \ref{sub-2.2}).

\br
The class of  minimally nice John metric spaces, which admit a roughly starlike Gromov hyperbolic quasihyperbolization, is very wide. For example, it includes (inner) uniform domains (more generally, uniform metric spaces),
simply connected John domains in the plane, and Gromov $\delta$-hyperbolic John domains in $\mathbb{R}^n$.
\er

\br
In view of the above, Theorem \ref{thm-1} states that all of the quasihyperbolic geodesics in the mentioned spaces are cone arcs. Moreover, Theorem \ref{thm-1} answers positively to question \ref{q-2} and also to question \ref{q-3} under a relatively mild condition.
\er

\br The main tool in the proof of Theorem \ref{thm-1} is the uniformization process of (Gromov) hyperbolic spaces, which was introduced by Bonk, Heinonen and Koskela in \cite{BHK}. They proved that each proper, geodesic and roughly starlike $\delta$-hyperbolic space is quasihyperbolically equivalent to a $c$-uniform space; see \cite[4.5 and 4.37]{BHK}. The uniformization process  of Bonk, Heinonen and Koskela  has many applications and is an important tool in many related papers, see e.g. \cite{BB03, KLM14}.
\er

From \cite[Theorem 3.22]{Vai05} it follows that every $\delta$-hyperbolic domain of $\IR^n$ is $K$-roughly starlike with $K$ depending only on $\delta$. Then we have the following corollary of Theorem \ref{thm-1}.

\bcor Every quasihyperbolic geodesic in an $a$-John, $\delta$-hyperbolic domain $D$ of $\IR^n$ is a $b$-cone arc with $b$ depending only on $a$ and $\delta$. \ecor

\br A proper domain $D$ in $\mathbb{R}^n$ is called $\delta$-{\it hyperbolic} for some $\delta\geq 0$, if $D$ has a Gromov hyperbolic quasihyperbolization. We remark that the above result is an improvement of \cite[Lemma $3.9$]{Guo15} whenever $\varphi(t)=Ct$ for some positive constant $C$. Also, we do not require the domain to be  bounded.
\er

\br There are many applications of the above mentioned classes of domains of $\mathbb{R}^n$ in the quasiconformal mappings and potential theory,
see e.g. \cite{BHK, CP17,GNV94,Guo15, NV}. A crucial ingredient in the related arguments is based on the fact that quasihyperbolic geodesics in Gromov hyperbolic John domains of $\mathbb{R}^n$ are inner uniform curves.
%
\er

As another motivation of this stude, we remark that Bonk, Heinonen and Koskela established the following characterization of Gromov hyperbolic domains on the $2$-sphere in \cite{BHK}.

\begin{Thm}\label{Thm-1} $($\cite[Theorem 1.12]{BHK}$)$ Gromov hyperbolic domains on the $2$-sphere are precisely the conformal images of inner uniform slit domains.
\end{Thm}

A {\it slit domain} is a proper subdomain $D$ of Riemann sphere such that each component of its complement is a point or a line segment parallel to the real or imaginary axis. It is well known that every domain in Riemann sphere is conformally equivalent to a slit domain. In \cite{BHK}, Bonk, Heinonen and Koskela also pointed out that their proof of Theorem \Ref{Thm-1} is ``surprisingly indirect, using among other things the theory of modulus and Loewner spaces as developed recently in \cite{HK}, plus techniques from harmonic analysis", and ask for an elementary proof as well.

In \cite{BB03}, Balogh and Buckley proved that a minimally nice metric space has a Gromov hyperbolic quasihyperbolization if and only if it satisfies the Gehring-Hayman condition and a ball separation condition. Their proof is also based on an analytic assumption that the space supports a suitable Poincar\'{e} inequality. Recently, Koskela, Lammi and Manojlovi\'{c} in \cite{KLM14} observed that Poincar\'{e} inequalities are not critical for this characterization of Gromov hyperbolicity, see \cite[Theorem 1.2]{KLM14}.

By using the above results, and as an application of Theorem \ref{thm-1}, we give the following simple geometric condition connecting the uniformity of a space to its other properties:

\begin{thm}\label{thm-2} Let $Q>1$ and let $(X,d,\mu)$ be a proper, $Q$-regular $A$-annularly quasiconvex length metric measure space. Let $D$ be a bounded proper subdomain of $X$. Then $D$ is uniform if and only if it is John or linearly locally connected, quasiconvex, and has a Gromov hyperbolic quasihyperbolization.
\end{thm}



\br With the aid of Theorem \ref{thm-1} and some auxiliary results obtained in \cite{KLM14}, the proof of Theorem \ref{thm-2} is essentially elementary and only needs the techniques from metric geometry and some estimates concerning the quasihyperbolic metrics. It is not difficult to find that Theorem \Ref{Thm-1} is a direct corollary of Theorem \ref{thm-2}.
\er

This paper is organized as follows. Section 2 contains notation and the basic definitions  and auxiliary lemmas. In Section 3, we will prove Theorem \ref{thm-1}. The proof of Theorem \ref{thm-2} is presented in Section 4.

\section{Preliminaries}

\subsection{Metric geometry} Let $(D, d)$ be a metric space, and let $B(x,r)$ and $\overline{B}(x,r)$ be the open ball and closed ball (of radius $r$ centered at the point $x$) in $D$, respectively. For a set $A$ in $D$, we use $\overline{A}$ to denote the metric completion of $A$ and $\partial A=\overline{A}\setminus A$ to be its metric boundary. A metric space $D$ is called {\it proper} if its closed balls are compact. Following terminology of \cite{BB03}, we call a locally compact, rectifiably connected noncomplete metric space $(D,d)$ {\it minimally nice}.

By a curve, we mean a continuous function $\gamma:$ $[a,b]\to D$. If $\gamma$ is an embedding of $I$,
it is also called an {\it arc}.
The image set $\gamma(I)$ of $\gamma$ is also denoted by $\gamma$. A curve $\gamma$ is called {\it rectifiably}, if the length $\ell_d(\gamma)<\infty$.
A metric space $(D, d)$ is called {\it rectifiably connected} if every pair of points in $D$ can be joined with a rectifiable
curve $\gamma$.

The length function associated with a rectifiable curve $\gamma$: $[a,b]\to D$ is $z_{\gamma}$: $[a,b]\to [0, \ell(\gamma)]$, given by
$z_{\gamma}(t)=\ell(\gamma|_{[a,t]})$.
For any rectifiable curve $\gamma:$ $[a,b]\to D$, there is a unique map $\gamma_s:$ $[0, \ell(\gamma)]\to D$ such that $\gamma=\gamma_s\circ z_{\gamma}$. Obviously,
$\ell(\gamma_s|_{[0,t]})=t$ for $t\in [0, \ell(\gamma)]$. The curve $\gamma_s$ is called the {\it arclength parametrization} of $\gamma$.

For a rectifiable curve $\gamma$ in $D$, the line integral over $\gamma$ of each Borel function $\varrho:$ $D\to [0, \infty)$ is
$$\int_{\gamma}\varrho ds=\int_{0}^{\ell(\gamma)}\varrho\circ \gamma_s(t) dt.$$

We say an arc $\gamma$ is {\it geodesic} joining $x$ and  $y$ in $D$ means that $\gamma$ is a map from an interval $I$ to $D$ such that $\gamma(0)=x$, $\gamma(l)=y$ and
$$\;\;\;\;\;\;\;\;d(\gamma(t),\gamma(t'))=|t-t'|\;\;\;\;\;\;\;\;\;\;\;\;\;\;\;\;\;\;\;\;\;\;\;\;\;\;\;\;\;\;\;\mbox{for all}\;\;t,t'\in I.$$

Every rectifiably connected metric space $(D, d)$ admits a natural (or intrinsic) metric, its  so-called length distance given by
$$\ell(x, y) := \inf\ell(\gamma)$$  where $ \gamma$   is a rectifiable curve joining $ x, y $ in $D.$
A metric space $(D, d)$ is a {\it length space} provided that $d(x, y) = \ell(x, y)$ for all points $x, y\in D$. It is also common to call such a $d$ an intrinsic distance function.

\subsection{Quasihyperbolic metric, quasigeodesics and solid arcs}\label{sub-2.2}
Suppose $\gamma $ is a rectifiable curve in a minimally nice space $(D,d)$, its {\it quasihyperbolic length} is the number:

$$\ell_{k_D}(\gamma)=\int_{\gamma}\frac{|dz|}{d_D(z)},
$$ where $d_D(z)=\dist(x,\partial D)$ is the distance from $z$ to the boundary of $D$.

For each pair of points $x$, $y$ in $D$, the {\it quasihyperbolic distance}
$k_D(x,y)$ between $x$ and $y$ is defined by
$$k_D(x,y)=\inf\ell_{k_D}(\gamma),
$$
where the infimum is taken over all rectifiable curves $\gamma$
joining $x$ to $y$ in $D$.
 We remark that the resulting space $(D,k_D)$ is complete, proper and geodesic (cf. \cite[Proposition $2.8$]{BHK}).

We recall the following basic estimates for quasihyperbolic distance that
first used by Gehring and Palka \cite[2.1]{GP76} (see also \cite[(2.3), (2.4)]{BHK}):
\be\label{li-1} k_D(x,y)\geq \log\Big(1+\frac{d(x,y)}
{\min\{d_D(x), d_D(y)\}}\Big)\geq \log|\frac{d_D(x)}{d_D(y)}|.\ee
In fact, more generally, we have
\be\label{li-2} \ell_{k_D}(\gamma)\geq \log\Big(1+\frac{\ell(\gamma)}
{\min\{d_D(x), d_D(y)\}}\Big)
\ee


Moreover, we have the following estimate:

\begin{lem}\label{newlemlabel}
Let $D$ be a minimally nice length space. Then for $x,y\in D$ with $d(x,y) < d_D(x)$, we have
$$k_D(x,y)\leq \frac{d(x,y)}{d_D(x)-d(x,y)}.$$
\end{lem}
\bpf  Let $0<\epsilon<\frac{1}{2}(d_D(x)-d(x,y))$. Since $D$ is a length space, there is a curve $\alpha$ joining $x$ and $y$ such that $\ell(\alpha)\leq d(x,y)+\epsilon$. Thus we have $\ell(\alpha)<d_D(x)$, which implies that $\alpha\subset B(x,d_D(x))\cap D$. Hence, we compute
$$k_D(x,y)\leq \int_{\alpha}\frac{|dz|}{d_D(z)}\leq \frac{\ell(\alpha)}{d_D(x)-\ell(\alpha)}<\frac{d(x,y)+\epsilon}{d_D(x)-d(x,y)-\epsilon}.$$
By letting $\epsilon\to 0$, we get the desired inequality.
\epf

\begin{defn} \label{def1.4}
 Suppose $\gamma$ is an arc in a minimally nice space $D$. The arc may be closed, open or half open. Let $\overline{x}=(x_0,$ $\ldots,$ $x_n)$,
$n\geq 1$, be a finite sequence of successive points of $\gamma$.
For $h\geq 0$, we say that $\overline{x}$ is {\it $h$-coarse} if
$k_D(x_{j-1}, x_j)\geq h$ for all $1\leq j\leq n$. Let $\Phi_{k_D}(\gamma,h)$
denote the family of all $h$-coarse sequences of $\gamma$. Set

$$z_{k_D}(\overline{x})=\sum^{n}_{j=1}k_D(x_{j-1}, x_j)$$ and
$$\ell_{k_D}(\gamma, h)=\sup \{z_{k_D}(\overline{x}): \overline{x}\in \Phi_{k_D}(\gamma,h)\}$$
with the agreement that $\ell_{k_D}(\gamma, h)=0$ if
$\Phi_{k_D}(\gamma,h)=\emptyset$. Then the number $\ell_{k_D}(\gamma, h)$ is the
{\it $h$-coarse quasihyperbolic length} of $\gamma$.  \end{defn}

\begin{defn} \label{def1.5} Let $D$ be  a minimally nice space. An arc $\gamma\subset D$
is {\it $(\nu, h)$-solid} with $\nu\geq 1$ and $h\geq 0$ if
$$\ell_{k_D}(\gamma[x,y], h)\leq \nu\;k_D(x,y)$$ for all $x$, $y\in \gamma$. \end{defn}

 Let $\lambda\geq 1$ and $\mu\geq 0$.
A curve $\gamma$ in $D$ is a {\it $(\lambda, \mu)$-quasigeodesic} if
$$\ell_{k_D}(x,y) \leq \lambda k_D(x,y)+\mu$$
for all $x,y\in \gamma.$ If $\lambda=1$, $\mu=0$, then $\gamma$ is a quasihyperbolic geodesic.

\begin{defn}Let $D$ and $D'$ be two minimally nice metric spaces. We say
that a homeomorphism $f: D\to D'$ is an {\it $M$-quasihyperbolic mapping},
or briefly  {\it $M$-QH}, if there exists a constant $M\geq 1$ such that for all $x$, $y\in D$,
$$\frac{1 }{M}k_D(x,y)\leq k_{D'}(f(x),f(y))\leq M\;k_D(x,y) .$$\end{defn}


In the following, we  use $x$, $y$, $z$, $\ldots$ to
denote the points in $D$, and $x'$, $y'$, $z'$, $\ldots$ the images
of $x$, $y$, $z$, $\ldots$ in $D'$, respectively, under $f$. For
arcs $\alpha$, $\beta$, $\gamma$, $\ldots$ in $D$, we also use
$\alpha'$, $\beta'$, $\gamma'$, $\ldots$ to denote their images in
$D'$.
Under  quasihyperbolic mappings, we have the following useful relationship between $(\lambda, \mu)$-quasigeodesics and solid arcs.

\begin{lem}\label{ll-001} Suppose that $G$ and $G'$ are minimally nice metric spaces. If $f:\;G\to G'$ is $M$-QH, and $\gamma$ is a $(\lambda, \mu)$-quasigeodesic in $G$, then there are constants $\nu=\nu(\lambda, \mu, M)$ and $h=h(\lambda, \mu, M)$ such that
 the image $\gamma'$ of $\gamma$ under $f$ is $(\nu,h)$-solid in $G'$.
\end{lem}
\bpf
Let $\gamma$ be a $(\lambda,\mu)$-quasigeodesic and let $$h=M(\lambda+\mu)\;\; \mbox{and}\;\; \nu=M^2(\lambda+\mu).$$
To show that $\gamma'$ is $(\nu,h)$-solid, we only need to verify that for $x$, $y\in \gamma$,
\be\label{new-eq-3}\ell_{k_{G'}}(\gamma'[x',y'],h)\leq\nu k_{G'}(x',y').\ee We prove this by considering two cases.
The first case is: $k_G(x,y)<1$. Then for $z$, $w\in\gamma[x, y]$, we have
$$k_{G'}(z',w')\leq Mk_G(z,w)\leq M(\lambda k_G(x,y)+\mu)<M(\lambda+\mu)=h,$$
and so
\be\label{ma-3}\ell_{k_{G'}}(\gamma'[x',y'],h)=0.\ee

Now, we consider the other case: $k_G(x,y)\geq 1$. Then with the aid of \cite[Theorem  4.9]{Vai6}, we have
\beq\label{ma-4}
\ell_{k_{G'}}(\gamma'[x',y'],h) &\leq&
 M\ell_{k_G}(\gamma[x,y])\leq M(\lambda k_G(x,y)+\mu)\\ \nonumber &\leq& M(\lambda+\mu)k_{G}(x,y) \leq  M^2(\lambda+\mu)k_{G'}(x',y').\eeq
It follows from \eqref{ma-3} and \eqref{ma-4} that \eqref{new-eq-3} holds, completing the proof.\epf


\subsection{Uniform spaces and John spaces}

In this subsection we first recall the definitions of John spaces, cone arcs and uniform spaces.  We also give some results related to some special  arcs which will be useful later in the proof of the main result.

\begin{defn} Let $a\geq 1$. A minimally nice space $(D,d)$ is called {\it $a$-John} if each pair of points $x,y\in D$ can be joined by a rectifiable arc $\alpha$ in $D$ such that for all $z\in \alpha$
$$\min\{\ell(\alpha[x,z]), \ell(\alpha[z,y])\}\leq a d_D(z),$$
where $\alpha[x,z]$ and $\alpha[z,y]$ denote two subarcs of $\alpha$ divided by the point $z$. The arc $\alpha$ is called an {\it $a$-cone} arc.
\end{defn}

\begin{defn} Let $c\geq 1$. A minimally nice space $(D,d)$ is called {\it $c$-uniform} if each pair of points $x,y\in D$ can be joined by a $c$-uniform arc. An arc $\alpha$ is called {\it $c$-uniform} if it is a  $c$-cone arc and satisfies the $c$-quasiconvexity, that is, $\ell(\alpha)\leq c d(x,y).$
\end{defn}

\begin{Lem}\label{Lem-uniform}$($\cite[(2.16)]{BHK}$)$\; If $D$ is a $c$-uniform metric space, then for all $x,y\in D$, we have
$$ k_{D}(x,y)\leq 4c^2 \log\Big(1+\frac{d(x,y)}{\min\{d_{D}(x),d_{D}(y)\}}\Big).$$\end{Lem}

The following properties of solid arcs in uniform metric spaces  is from \cite{LVZ2} which will be used in our proofs.

\begin{Lem}\label{Lem13''}$($\cite[Lemma 3]{LVZ2}$)$\, Suppose that $D$  is a $c$-uniform space, and that $\gamma$ is a $(\nu,h)$-solid arc in $D$ with endpoints $x$, $y$. Let $d_D(x_0)=\max_{p\in \gamma}d_D(p)$. Then there exist constants
$a_1=a_1( c, \nu, h)\geq 1$ and $a_2=a_2(c, \nu, h)\geq 1$ such that
\begin{enumerate} \item
$\diam(\gamma[x,u])\leq a_1 d_D(u)$ for $u\in \gamma[x,x_0],$ and $\diam(\gamma[y,v])\leq a_1 d_D(v)$ for $v\in \gamma[y, x_0]$;
\item  $\diam(\gamma)\leq \max\big\{a_2 d(x,y), 2(e^h-1)\min\{d_D(x),d_D(y)\}\big\}.$
\end{enumerate}
 \end{Lem}

%

Next we discuss the properties of cone arcs. 

\begin{lem}\label{eq-8} Let $\alpha[x,y]$ be an $a$-cone arc in $D$ and let $z_0$ bisect the arclength of
$\alpha[x,y]$. Then for each $z_1$, $z_2\in\alpha[x,z_0]$ $($or
$\alpha[y,z_0]$$)$ with $z_2\in \alpha[z_1,z_0]$, we have $$k_D(z_1,z_2)\leq  \ell_k(\alpha[z_1,z_2])\leq 2a \log\big(1+\frac{2\ell(\alpha[z_1,z_2])}{d_D(z_1)}\big)$$ and $$\ell_{k}(\alpha[z_1,z_2])\leq
4a^2k_{D}(z_1,z_2)+4a^2.$$
\end{lem}

\bpf By symmetry, we only need to verify the assertion in the case $z_1$, $z_2\in\alpha[x,z_0]$. To this end, for $z_2\in \alpha[z_1,z_0]$ be given, we have
$$d_D(z_2)\geq \frac{1}{a}\ell(\alpha[z_1,z_2]).$$

If $z_2\subset
B(z_1, \frac{1}{2}d_D(z_1))$, thus one  finds that $d_D(z_2)\geq
\frac{1}{2}d_D(z_1)$.
Otherwise, we have
$d_D(z_2)\geq\frac{1}{2a}d_D(z_1)$. Hence in both cases we obtain
$$d_D(z_2)\geq \frac{1}{4a}[2\ell(\alpha[z_1,z_2])+d_D(z_1)],$$
 which yields that
\begin{eqnarray*}k_{D}(z_1,z_2) &\leq& \ell_k(\alpha[z_1,z_2])= \int_{\alpha[z_1,z_2]}\frac{|dz|}{d_D(z)}\\ \nonumber
&\leq& 2a\log\Big(1+\frac{2\ell(\alpha[z_1,z_2])}{d_D(z_1)}\Big)\\ \nonumber
&\leq& 4a^2\log\Big(1+\frac{d_D(z_2)}{d_D(z_1)}\Big)\\ \nonumber &\leq&
4a^2k_{D}(z_1,z_2)+4a^2,\end{eqnarray*}
as desired.
\epf

\blem \label{lem13-0-0}Suppose that
 $f: D\to D'$ is an $M$-QH from an $a$-John minimally nice space $D$ to a $c$-uniform space $D'$. Let $\alpha$ be an $a$-cone arc in $D$ with end points $x$ and $y$,  $z_0$ bisect the arclength of
$\alpha$, and let $d_{D'}(v'_1)=\max\{d_{D'}(u'): u'\in\alpha'[x',z'_0]\}$ and
$d_{D'}(v'_2)=\max\{d_{D'}(u'): u'\in\alpha'[y',z'_0]\}$. Then there is a constant $a_3=a_3(a,c,M)$ such that
 \begin{enumerate}

\item  for each  $z'\in \alpha'[x', v'_1]$,
$d' (x',z')\leq a_3\;d_{D'}(z')$
 and for each  $z'\in \alpha'[v'_1, z'_0]$, $d' (z'_0,z')\leq a_3\;d_{D'}(z')$.
 \item  for each  $z'\in \alpha'[y', v'_2]$,
$d' (y',z')\leq a_3\;d_{D'}(z')$
 and for each  $z'\in \alpha'[v'_2, z'_0]$, $d' (z'_0,z')\leq a_3\;d_{D'}(z')$. \end{enumerate}
\elem
\bpf  First, in the light of Lemma \ref{eq-8}, we see that $\alpha[x,z_0]$ and $\alpha[z_0,y]$ are $(4a^2,4a^2)$-quasigeodesics. Since $f: D\to D'$ is $M$-QH, we thus know from Lemma \ref{ll-001} that $\alpha'[x',z'_0]$ and $\alpha'[z'_0,y']$ are solid arcs. Moreover, by the choices of $v_1'$ and $v_2'$, $(1)$ and $(2)$ follows from Lemma \Ref{Lem13''}.
\epf

\subsection{Uniformization theory of Bonk, Heinonen and Koskela }
Let $(X,d)$ be a geodesic metric space and let $\delta\geq 0$. If for all triples of geodesics $[x,y], [y,z], [z,x]$ in $(X,d)$ satisfies: every point in $[x,y]$ is within distance $\delta$ from $[y,z]\cup [z,x]$, then the space $(X,d)$ is called a {\it $\delta$-hyperbolic space}. For simplicity, in the rest of this paper when we say that a minimally nice space $X$ is {\it Gromov
hyperbolic} we mean that the space is $\delta$-hyperbolic with respect to the quasihyperbolic metric for some nonnegative constant $\delta$.

In \cite{BHK}, Bonk, Heinonen and Koskela  introduced the concept of rough starlikeness of a Gromov hyperbolic space with respect to a given base point. Let $X$ be a proper, geodesic $\delta$-hyperbolic space, and let $w\in X$, we say that $X$ is {\it $K$-roughly starlike} with respect to $w$ if for each $x\in X$ there is some point $\xi\in\partial^* X$ and a geodesic ray $\gamma=[w,\xi]$ with $\dist(x,\gamma)\leq K$.

They also proved that both bounded uniform spaces and every hyperbolic domain $($a domain equipped with its quasi-hyperbolic metric is a Gromov hyperbolic space$)$ in $\IR^n$ are roughly starlike. It turns out that this property serves as an important tool in several research, for instance \cite{BB03}, \cite{ZR} and \cite{KLM14}.

Next we recall the conformal deformations which were introduced by Bonk, Heinonen and Koskela (cf. \cite[Chapter $4$]{BHK}). Let $(X,d)$ be a minimally nice space and $w\in X$. Consider the family of conformal deformations of $(X,k)$ by the densities
$$\rho_\epsilon(x)=e^{-\epsilon k(x,w)}\;\;(\epsilon>0).$$
For $u$, $v\in X$, let
$$d_\epsilon(u,v)=\inf\int_{\gamma} \rho_\epsilon ds_k,$$
where $ds_k$ is the arc-length element with respect to the metric $k$ and the infimum is taken over all rectifiable curves $\gamma$ in $X$ with endpoints $u$ and $v$.
Then $d_\epsilon$ are metrics on $X$, and we denote the resulting metric spaces
by $X_\epsilon=(X,d_\epsilon)$.

The next result shows that the deformations $X_{\epsilon}$ are uniform spaces and  each proper, geodesic and roughly starlike $\delta$-hyperbolic space is {\it quasihyperbolically equivalent} to a $c$-uniform space; see \cite[Propositions 4.5 and 4.37]{BHK}.

\begin{Lem}\label{lem-1}$($\cite[Propositions $4.5$ and $4.37$]{BHK} or \cite[Lemma $4.12$]{BB03}$)$
Suppose $(X,d)$ is minimally nice, locally compact and that $(X,k)$ is both $\delta$-Gromov hyperbolic and $K$-roughly starlike, for some $\delta\geq 0$, $K>0$. Then $X_\epsilon$ has diameter at most $2/\epsilon$ and there are positive numbers $c, \epsilon_0$  depending only on $\delta, K$ such that $X_\epsilon$ is $c$-uniform for all $0<\epsilon\leq \epsilon_0$. Furthermore, there exists $c_0=c_0(\delta,K)\in(0,1)$ such that the quasihyperbolic metrics $k$ and $k_\epsilon$ satisfy the quasi-isometric condition
$$c_0\epsilon k(x,y)\leq k_\epsilon(x,y)\leq e \epsilon k(x,y).$$
\end{Lem}

\bigskip

\section{The proof of Theorem \ref{thm-1}}

Let $(D,d)$ be a minimally nice $a$-John metric space and $(D,k)$ $K$-roughly starlike, $\delta$-hyperbolic where $k$ is the quasihyperbolic metric of $D$. Then by Lemma \Ref{lem-1}, we know that there is a positive number $\epsilon=\epsilon(\delta)$ such that $(D,d_{\epsilon})$ is a $c$-uniform metric space and the identity map from $(D,d)$ to $(D,d_{\epsilon})$ is $M$-QH, where $c$ and $M$ depend only on $\delta$ and $K$. For simplicity, we denote $D=(D,d)$, $(D',d')=(D,d_{\epsilon})$ and $f$ the identity map from $D$ to $D'$.

We may assume without loss of generality that $D$ is a length space, because the length of an arc and the quasihyperbolic metrics associated to the original metric and the length metric coincide.

%

 Fix $z_1$, $z_2\in D$ and let $\gamma$ be a quasihyperbolic geodesic joining
$z_1$, $z_2$ in $D$. Let $b=4a_4e^{a_4}$, $a_4=a_5^{8c^2M}$,
$a_5=a_6^{4a^2M}$ and $a_6=(8a_1^2a_3)^{16c^2M}a^2$, where $a_1$ and $a_3$ are the constants from Lemmas \Ref{Lem13''} and \ref{lem13-0-0}, respectively. In the following, we shall prove that $\gamma$ is a
$b$-cone arc, that is, for each $y\in\gamma$,
$$\min\{\ell(\gamma[z_1, y]),\; \ell(\gamma[z_2,
y])\}\leq b\,d_D(y).$$ 

Let $x_0\in \gamma$ be a point such that
$d_D(x_0)=\max\limits_{z\in \gamma}d_D(z).
$
By symmetry, we only need to prove that for $y\in\gamma[z_1,x_0]$, \beq \label{John}\ell(\gamma[z_1, y])\leq b\,d_D(y).\eeq

To this end, let $m \geq 0$
be an integer such that
$$2^{m}\, d_D(z_1)
\leq d_D(x_0)< 2^{m+1}\, d_D(z_1).
$$
And let $y_0$ be the first point in $\gamma[z_1,x_0]$ from $z_1$ to
$x_0$ with
$$d_D(y_0)=2^{m}\, d_D(z_1).
$$
Observe that if $d_D(x_0)=d_D(z_1)$, then $y_0=z_1=x_0$.

Let $y_1=z_1$. If $z_1=y_0$, we let $y_2=x_0$. It is possible that
$y_2=y_1$. If $z_1\not= y_0$, then we let $y_2,\ldots ,y_{m+1}$ be
the points such that for each $i\in \{2,\ldots,m+1\}$, $y_i$
denotes the first point in $\gamma[z_1,x_0]$ from $y_1$ to $x_0$
satisfying
$$d_D(y_i)=2^{i-1}\, d_D(y_1).$$
Then $y_{m+1}=y_0$. We let $y_{m+2}=x_0$. It is possible that
$y_{m+2}=y_{m+1}=x_0=y_0$. This possibility occurs once
$x_0=y_0$.

From the choice of $y_i$ we observe that for $y\in \gamma[y_i,y_{i+1}]$ $(i\in\{1, 2, \ldots, m+1\})$, \be\label{li-newadd-1} d_D(y)<d_D(y_{i+1})=2d_D(y_i)\ee
and so for all $i\in\{1, 2, \ldots, m+1\}$,
\be\label{li-newadd-2} k_{D}(y_i,y_{i+1}) =\ell_k(\gamma[y_i,y_{i+1}])\geq \frac{\ell(\gamma[y_i,y_{i+1}])}{2d_D(y_i)}.\ee

To prove Theorem \ref{thm-1}, we shall estimate upper bound of the quasihyperbolic distance between $y_i$ and $y_{i+1}$, which state as follows.

\begin{lem}\label{eq-0}For each $i\in \{1,\ldots, m+1\}$, $k_{D}(y_i,y_{i+1})\leq a_4$.\end{lem}

We note that Theorem \ref{thm-1} can be obtained from Lemma \ref{eq-0} as follows.
First, we observe from \eqref{li-newadd-2} and Lemma \ref{eq-0} that for all $i\in\{1,\ldots, m+1\}$,
\be\label{li-1'} \ell(\gamma[y_i,y_{i+1}])\leq
 2a_4 \,d_D(y_i).\ee

Further, for each $y\in \gamma[y_1,x_{0}]$, there is some $i\in
\{1,\ldots,m+1\}$ such that $y\in \gamma[y_i,y_{i+1}]$.  It follows from \eqref{li-1} that
$$ \log \frac{d_D(y_i)}{d_D(y)}\leq k_D(y,y_i)\leq
\, k_D(y_i,y_{i+1})\leq a_4 ,$$ whence
$$d_D(y_i)\leq e^{ a_4 }d_D(y).$$
From which and (\ref{li-1'}) it follows that
\beq\label{eq(li-3)} \ell(\gamma[z_1,y])&=&
\ell(\gamma[y_1,y_2])+\ell(\gamma[y_2,y_3])+\ldots+\ell(\gamma[y_i,y])
\\
\nonumber &\leq& 2a_4 (d_D(y_1)+d_D(y_2)+\ldots+d_D(y_i))\\ \nonumber
&\leq& 4a_4 \,d_D(y_i)\leq 4a_4 e^{a_4 }\,d_D(y),\eeq as desired. This proves \eqref{John} and so Theorem \ref{thm-1} follows.

Hence to complete the proof of Theorem \ref{thm-1}, we only need to prove Lemma \ref{eq-0}.

\subsection{The proof of Lemma \ref{eq-0}}
Without loss of generality, we may assume that
$d_{D'}(y'_i)\leq d_{D'}(y'_{i+1})$. We note that if  $ d (y_i, y_{i+1})<\frac{1}{2}d_D(y_i),$ then by Lemma \ref{newlemlabel} we have
$$k_D(y_i, y_{i+1})\leq 1,$$
as desired. Therefore, we assume in the following that
\beq\label{eq(4-2)}d (y_i, y_{i+1})\geq \frac{1}{2}d_D(y_i).\eeq

Let $\alpha_i$ be an
$a$-cone arc joining $y_i$ and $y_{i+1}$ in $D$ and let $v_i$ bisect the arclength of
$\alpha_i$. Then Lemma \ref{eq-8} implies that
\beq\label{hl-eq(4-1-2)}\;\;\;\;\;k_{D}(y_i,y_{i+1})&\leq&
k_{D}(y_i,v_i)+k_{D}(v_i,y_{i+1})\\ \nonumber &\leq& 2a\bigg(\log \Big( 1+\frac{2\ell(\alpha_i[y_i,v_i])}
{d_D(y_i)}\Big)+\log
\Big( 1+\frac{2\ell(\alpha_i[y_{i+1},v_i])} {d_D(y_{i+1})}\Big)\bigg)\\ \nonumber &\leq& 4a\log \Big(
1+\frac{\ell(\alpha_i)} {d_D(y_i)}\Big). \nonumber
\eeq


Now we divide the proof of Lemma \ref{eq-0} into two cases.

\bca \label{ca1}  $\ell(\alpha_i)< a_5  d (y_i, y_{i+1}).$\eca
Then by \eqref{li-newadd-2} and \eqref{hl-eq(4-1-2)} we compute
\beq\label{eq(h-h-4-2')}
 \frac{d(y_i,y_{i+1})}{2d_D(y_i)}&\leq& k_{D}(y_i,y_{i+1})
  \leq 4a \log \Big( 1+\frac{\ell(\alpha_i)} {d_D(y_i)}\Big) \\ \nonumber
  &\leq&  4a \log \Big( 1+\frac{ a_5d(y_i,y_{i+1})} {d_D(y_i)}\Big).\eeq A necessary condition for \eqref{eq(h-h-4-2')} is
 $$ d (y_i,y_{i+1})\leq a_5^2\,d_D(y_i).$$
Hence we deduce from (\ref{eq(h-h-4-2')}) that $k_{D}(y_i,y_{i+1})\leq
a_4$, as desired.

\bca \label{ca2}  $\ell(\alpha_i)\geq a_5  d (y_i, y_{i+1}).$\eca

We prove in this case by contradiction. Suppose on the contrary that
\beq\label{eq(h-4-2)}k_{D}(y_i,y_{i+1})> a_4.\eeq
Then by Lemma \Ref{Lem-uniform}, we get
\begin{eqnarray*}a_4<k_{D}(y_i,y_{i+1})\leq M k_{D'}(y'_i,y'_{i+1})
\leq 4c^2M\log\Big(1+\frac{d' (y'_i,y'_{i+1})}{d_{D'}(y'_i)}\Big),\end{eqnarray*}
and so
\beq\label{eq(h-4-1')}d' (y'_i,y'_{i+1})\geq a_5d_{D'}(y'_i).\eeq
Therefore, by the choice of $v_i\in\alpha_i$ we obtain
$$d_D(v_i)\geq \frac{\ell(\alpha_i)}{2a}\geq \frac{a_5}{2a} d (y_i,y_{i+1})>a_6\, d (y_i,y_{i+1}),$$
we deduce from which and \eqref{eq(4-2)} that there exists a point
$v_{i,0}\in \alpha_i[y_i,v_i]$ such that
\be\label{eq-11} d_D(v_{i,0})=a_6\, d (y_i,y_{i+1}).\ee
Moreover, we claim that
\be\label{claim1}k_{D}(y_i,v_{i,0})\leq \frac{1}{a_5}k_{D}(y_i,y_{i+1}).\ee
Otherwise, we would see from  Lemma \ref{eq-8} and \eqref{eq-11} that
 \begin{eqnarray*}k_{D}(y_i,y_{i+1})&<&  a_5  k_{D}(y_i,v_{i,0})\leq 4aa_5 \log\Big(1+\frac{\ell(\alpha_i[y_{i},v_{i,0}])}{d_D(y_i)}\Big)
\\ \nonumber
&\leq& 4aa_5 \log\Big(1+\frac{ad(v_{i,0})}{d_D(y_i)}\Big)
\leq
4a^2a_5a_6\log\Big(1+\frac{ d (y_i,y_{i+1})}{d_D(y_i)}\Big),
\end{eqnarray*} which together with \eqref{li-newadd-2} show that $$\frac{d(y_i,y_{i+1})}{d_D(y_i)}\leq 8a^2a_5a_6\log\Big(1+\frac{ d (y_i,y_{i+1})}{d_D(y_i)}\Big).$$
A necessary condition for the above inequality is   $$ d (y_i,y_{i+1})\leq a_5^2\,d_D(y_i).$$
This shows that $k_{D}(y_i,y_{i+1})\leq
a_4$, which contradicts $\eqref{eq(h-4-2)}$. Thus we get (\ref{claim1}).

Then it follows from Lemma \Ref{Lem-uniform}, and \eqref{claim1} that
\begin{eqnarray*} k_{D'}(y'_i,v'_{i,0})&<& Mk_{D}(y_i,v_{i,0})
  \leq\frac{M}{a_5}k_{D}(y_i,y_{i+1}) \\  &\leq& \frac{M^2}{ a_5}k_{D'}(y'_i,y'_{i+1})
\leq \frac{4c^2M^2}{ a_5}\log\Big(1+\frac{d' (y'_i,y'_{i+1})}{d_{D'}(y'_i)}\Big).
\end{eqnarray*}
Hence, by using an elementary compute we see from \eqref{li-1} and \eqref{eq(h-4-1')} that \begin{eqnarray*}
\log \Big(1+\frac{d' (y'_i,v'_{i,0})}{d_{D'}(y'_i)}\Big) \leq k_{D'}(y'_i,v'_{i,0})\leq \log\Big(1+\frac{d' (y'_i,y'_{i+1})}{a_5d_{D'}(y'_i)}\Big),
\end{eqnarray*}
which implies that
\beq\label{eq(hl-41-5)}d' (y'_i,v'_{i,0})<
\frac{1}{a_5}d' (y'_i,y'_{i+1}).\eeq
Moreover, we deduce from (\ref{eq(hl-41-5)}) and
(\ref{eq(h-4-1')}) that
\beq\label{eq--2} d_{D'}(v'_{i,0})\leq
d' (y'_i,v'_{i,0})+d_{D'}(y'_i)\leq
\frac{2}{a_5}d' (y'_i,y'_{i+1}).\eeq

We recall that $v_i$ is the point in the cone arc $\alpha_i[y_i,y_{i+1}]$ which bisect the arclength of $\alpha_i$. Next we need to estimate the location  of the image point $v'_i$ in $\alpha'_i$. We claim that

\begin{cl}\label{eq--6}
$d'(y'_i,v'_i)<\frac{1}{2}d' (y'_i,y'_{i+1}).$ \end{cl}

We prove this claim by a method of contradiction. Suppose on the contrary that \be\label{neweqlabel}d' (y'_i,v'_i)\geq
\frac{1}{2}d' (y'_i,y'_{i+1}).\ee Let $u'_{0,i}\in\gamma'[y'_{i},
y'_{i+1}]$ be a point satisfying $$d_{D'}(u'_{0,i})=\max\{d_{D'}(w'):w'\in\gamma'[y'_{i},
y'_{i+1}]\}.$$
Then we see from Lemma \Ref{Lem13''}
that
\be\label{e---1} d_{D'}(u'_{0,i})\geq \frac{1}{a_1}\max\{d' (y'_{i+1},u'_{0,i}),
d' (u'_{0,i},y'_i)\} \geq
\frac{d' (y'_i,y'_{i+1})}{2a_1}.\ee This together with (\ref{eq(h-4-1')}) shows that there  exists some point
$y'_{0,i}\in \gamma'[y'_i,u'_{0,i}]$ satisfying
\beq\label{eq(W-l-6-1)}d_{D'}(y'_{0,i})=\frac{d' (y'_i,y'_{i+1})}{2a_1}.
\eeq
It follows from Lemma  \Ref{Lem13''} that\beq\label{eq(W-l-6-1add)}d' (y'_i,y'_{0,i})\leq
a_1\,d_{D'}(y'_{0,i}).\eeq

%

Let $v'_0\in\alpha'_i[y'_{i}, v'_{i}]$ satisfy $d_{D'}(v'_0)=\max\{d_{D'}(u'):u'\in\alpha'_i[y'_{i}, v'_{i}]\}$, see Figure \ref{fig01}. Then we see from Lemma
\ref{lem13-0-0} that for each $z'\in \alpha'_i[ v'_i, v'_0]$, \beq\label{cla-3}d' (v'_i,z')\leq a_3 d_{D'}(z').\eeq On the other hand, we recall that $v'_{i,0}$ is the point such that $v_{i,0}\in \alpha_i[y_i,v_i]$ and  satisfying \eqref{eq-11} and \eqref{eq(hl-41-5)}. Then by \eqref{eq(hl-41-5)} and \eqref{eq--2} we have
\begin{eqnarray*}d' (v'_i,v'_{i,0})&\geq&
d' (v'_i,y'_i)-d' (v'_{i,0},y'_i)\geq (\frac{1}{2}-\frac{1}{a_5})d' (y'_i,y'_{i+1})>a_3 d_{D'}(v'_{i,0}).
\end{eqnarray*}
That means $v'_0\in \alpha'_i[ v'_{i,0}, v'_i]$.

Moreover, we know from Lemma
\ref{lem13-0-0} and \eqref{neweqlabel} that $$d_{D'}(v'_0)\geq \frac{1}{a_3}\max\{d' (v'_{i},v'_0),
d' (v'_0,y'_i)\}\geq \frac{d' (y'_i,v_i')}{2a_3}\geq
\frac{d' (y'_i,y'_{i+1})}{4a_3}.$$
Hence, it follows from  (\ref{eq--2})  that there exists some point $u'_0\in
\alpha'_i[v'_{i,0},v'_{0}]$ such that \beq\label{eq(W-l-6-2)}
d_{D'}(u'_0)=\frac{d' (y'_i,y'_{i+1})}{4a_3},\eeq and so Lemma \ref{lem13-0-0} leads to
$$d' (y'_i,u'_0)\leq a_3\,d_{D'}(u'_0).$$
This together with \eqref{eq(W-l-6-1)},  \eqref{eq(W-l-6-1add)} and \eqref{eq(W-l-6-2)} show  that
$$d' (u'_0,y'_{0,i})\leq d' (u'_0,y'_i)+d' (y'_i,y'_{0,i})\leq 3a_3d_{D'}(u'_0).$$
Now we are ready to finish the proof of Claim \ref{eq--6}.  It follows from \eqref{li-1} and Lemma \Ref{Lem-uniform} that
\begin{eqnarray*}
\log \frac{d_D(u_0)}{d_D(y_{0,i})}&\leq& k_{D}(y_{0,i},u_0)
\leq M k_{D'}(y'_{0,i},u'_0) \\ \nonumber
&\leq& 4c^2M\log\Big(1+\frac{d' (u'_0,y'_{0,i})}{\min\{d_{D'}(u'_0),
d_{D'}(y'_{0,i})\}}\Big)\\ \nonumber &<&4c^2M\log (1+3a_3),
\end{eqnarray*}
which yields that
\beq\label{eq(W-l-6-4)}d_D(u_0)\leq (1+3a_3)^{4c^2M}d_D(y_{0,i})<a_6d_D(y_{0,i}).\eeq

On the other hand, by Lemma \ref{eq-8} we can get
 \begin{eqnarray*}
  k_{D}(v_{i,0},u_0)\geq 4{a}^2\log\Big(1+\frac{d_D(u_0)}{d_D(v_{i,0})}\Big)\end{eqnarray*}

and by \eqref{li-1}, \eqref{eq--2} and \eqref{eq(W-l-6-2)} we have
that
\begin{eqnarray*}  k_{D}(v_{i,0},u_0)  \geq k_{D'}(v'_{i,0},u'_0)
\geq\log\frac{d_{D'}(u'_0)}{d_{D'}(v'_{i,0})}
\geq\log\frac{a_5}{8a_3},\end{eqnarray*} which yields $$d_D(u_0)\geq
a_6d_D(v_{i,0}).$$ Therefore, we infer from \eqref{eq(4-2)} and
\eqref{eq-11} that
\begin{align*} d_D(u_0)\geq  a_6d_D(v_{i,0})=a_6^2 d (y_i, y_{i+1})
\geq \frac{ a_6^2}{4}d_D(y_{i+1})\geq \frac{
a_6^2}{4}d_D(y_{0,i}),\end{align*} which contradicts
(\ref{eq(W-l-6-4)}). Hence Claim \ref{eq--6} holds.

\medskip

Now we continue the proof of Lemma \ref{eq-0}.
 We first see from Claim \ref{eq--6} that
$$d' (y'_{i+1},v'_i)\geq d'(y'_i,y'_{i+1})-d' (y'_i,v'_i)>
\frac{d' (y'_i,y'_{i+1})}{2}\geq d' (y'_i,v'_i).$$ Let $q'_0\in
\alpha'_i[y'_i,v'_i]$ and $u'_1\in\alpha'_i[y'_{i+1},v'_i]$ be points such that
\beq\label{112}\;\;\;\;\frac{d' (y'_i,v'_i)}{2a_3}= d' (q'_0,v_i')\,\,\;{\rm and}\,\,\; \frac{d'(y'_i,v'_i)}{2a_3}= d' (u'_1,v_i').\eeq
Then \begin{align*} d'(y'_i,q_0')\geq d'(y'_i,v_i')-d'(q'_0,v_i')=(2a_3-1)d'(q'_0,v_i')>d'(q'_0,v_i')\end{align*}
and \begin{align*} d'(y'_{i+1},u_1')>d'(u'_1,v_i').\end{align*}
Thus we get from  Lemma \ref{lem13-0-0} that
\beq\label{eq(W-l-6'-0)}\;\;\;\;\;\;\;\;\;\;\;\; d_{D'}(q'_0)\geq \frac{d' (q'_0,v_i')}{a_3}\geq
\frac{d' (y'_i,v'_i)}{2a^2_3} \,\; \mbox{and}\;\,
d_{D'}(u'_1)\geq\frac{d' (u'_1,v_i')}{a_3}\geq
\frac{d' (y'_i,v'_i)}{2a_3^2}.\eeq
Then it follows from Lemma \Ref{Lem-uniform}, \eqref{li-1}, \eqref{112} and \eqref{eq(W-l-6'-0)} that
\beq \label{eq--7}
\Big|\log\frac{d_D(u_1)}{d_D(q_0)}\Big| &\leq&k_{D}(u_1,q_0) \leq M k_{D'}(u'_1, q'_0)
\\\nonumber &\leq& 4c^2M\log\Big(1+\frac{d' (u'_1,q'_0)}{\min\{d_{D'}(q'_0), d_{D'}(u'_1)\}}\Big)
\\\nonumber &\leq& 4c^2M\log\Big(1+\frac{d' (u'_1,v'_i)+d' (v'_i,q'_0)}{\min\{d_{D'}(q'_0), d_{D'}(u'_1)\}}\Big)
\\\nonumber &\leq& 4c^2M\log (1+2a_3 ),
\eeq
which implies that
\beq\label{eq(W-l-6'-2)}\frac{d_D(u_1)}{(1+2a_3 )^{4c^2M}}\leq
d_D(q_0)\leq (1+2a_3 )^{4c^2M}e^Cd_D(u_1).\eeq

%

On the other hand, by \eqref{eq(h-4-1')}, \eqref{e---1} and Claim \ref{eq--6} we have
$$d' (u'_{0,i},y'_i)\geq d_{D'}(u'_{0,i})-d_{D'}(y'_i)\geq (\frac{1}{2a_1}-\frac{1}{a_5})d' (y'_{i+1},y'_i)>\frac{1}{2a_1}d' (y'_{i},v'_i).$$
Then there exists $p'_0\in
\gamma'[y'_i,u'_{0,i}]$ such that
\beq\label{132}d' ( y'_i,p'_0)=\frac{d' (y'_i,v'_i)}{2a_1},\eeq
see Figure \ref{fig02}.
This combined with \eqref{112} and Lemma \Ref{Lem13''} shows that $$d' (p'_0,q'_0)\leq  d' ( p'_0,y'_i)+d' (y'_i,v'_i)+d' (v'_i,q'_0)\leq (1+\frac{1}{a_1}+\frac{1}{a_3})d' (y'_i,v'_i),$$
and
$$d' (y'_i,p'_0 )\leq a_1 d_{D'}(p'_0).$$ Then \eqref{eq(W-l-6'-0)} and \eqref{132}  we have $$\min\{d_{D'}(q'_0),d_{D'}(p'_0)\}\geq \min\{\frac{1}{2a_1^2},\frac{1}{2a_3^2}\}d'(y_i',v_i')
>\frac{1}{2a_1^2a_3^2}d'(y_i',v_i').$$ Therefore, Lemma \ref{newlemlabel} and \eqref{li-1} lead to
\begin{eqnarray*}\log \frac{d_D(q_0)}{d_D(p_{0})}&\leq& k_{D}(q_0, p_{0})
\leq M k_{D'}(q'_0,p'_0)\\
\nonumber &\leq&
4c^2M\log\Big(1+\frac{d' (p'_0,q'_0)}{\min\{d_{D'}(q'_0),d_{D'}(p'_0)\}}\Big)
\\
\nonumber &\leq& 4c^2M\log(6a_1^2a_3^2).\end{eqnarray*}
 We infer
from   \eqref{eq-11} that
\begin{eqnarray}\label{eq-new-add2}d_D(q_0)&\leq&
(6a_1^2a_3^2)^{4c^2M} d_D(p_0)\\ \nonumber&\leq&
2(6a_1^2a_3^2)^{4c^2M} d_D(y_i)
\\ \nonumber&\leq& 2(6a_1^2a_3^2)^{4c^2M} d (y_i,y_{i+1}).\end{eqnarray}

Finally, it follows from  Lemma \ref{eq-8} and the choice of $q_0$ and $u_1$ that $$k_{D}(y_i, q_0)\leq 4a^2\log\Big(1+\frac{d_D(
q_0)}{d_D(y_i)}\Big)\;\;{\rm
and }\;\;k_{D}(u_i, y_{i+1})\leq 4a^2\log\Big(1+\frac{d_D(u_i)}{d_D(y_{i+1})}\Big).$$
Then by Lemma \Ref{Lem-uniform}, \eqref{eq--7}, \eqref{eq(W-l-6'-2)} and \eqref{eq-new-add2} we get
\beq\label{eq(W-l-6'-2')}
  k_{D}(y_i, y_{i+1})&\leq&   k_{D}(y_i, q_0)+k_{D}(q_0, u_1)+k_{D}(u_1, y_{i+1})
\\ \nonumber &\leq& 4a^2 \log\Big(1+\frac{d_D(
q_0)}{d_D(y_i)}\Big)+4A^2M \log\Big(1+2a_3\Big)
\\ \nonumber &&+4a^2 \log\Big(1+\frac{d_D(
u_1)}{d_D(y_{i+1})}\Big)\\
\nonumber &<& a_5 \log\Big(1+\frac{ d (y_i,
y_{i+1})}{d_D(y_i)}\Big),\eeq which together with \eqref{li-newadd-2} show that $$\frac{d (y_i,y_{i+1})}{2d_D(y_i)}\leq a_5 \log\Big(1+\frac{ d (y_i,y_{i+1})}{d_D(y_i)}\Big).$$ A necessary condition for
this inequality is $ d (y_i, y_{i+1})\leq a_5^2d_D(y_i)$.
Hence by (\ref{eq(W-l-6'-2')}), we know that
$$k_{D}(y_i, y_{i+1})\leq a_5\log(1+a_5^2)<a_4,$$
 which   contradicts \eqref{eq(h-4-2)}.
Therefore,  we obtain Lemma \ref{eq-0} and so Theorem \ref{thm-1}.

\bigskip

\section{The proof of Theorem \ref{thm-2}}

In this section, we will prove Theorem \ref{thm-2} by means of Theorem \ref{thm-1} and some results demonstrated in \cite{KLM14}. We begin by recalling necessary definitions and results.

\bdefe Let $(X, d,\mu)$ be a metric measure space. Given $Q> 1$, we say that
$X$ is {\it $Q$-regular} if there exists a constant $C>0$ such that
for each $x\in X$ and $0<r\leq \diam(X)$,
$$C^{-1}r^Q\leq \mu(B(x,r))\leq Cr^Q.$$\edefe

\bdefe
Let $(X,d)$ be a locally compact and rectifiably connected metric space, $D\subset X$ be a domain (an open rectifiably connected set), and $C_{gh}\geq 1$ be a constant. We say that $D$ satisfies the {\it $C_{gh}$-Gehring-Hayman inequality}, if for all $x$, $y$ in $D$ and for each
quasihyperbolic geodesic $\gamma$ joining $x$ and $y$, we have
$$\ell(\gamma)\leq C_{gh}\ell(\beta_{x,y}),$$
where $\beta_{x,y}$ is any other curve joining $x$ and $y$ in $D$. In other words,  quasihyperbolic geodesics are essentially the shortest curves in $D$.
\edefe

\bdefe
Let $(X,d)$ be a metric space, $D\subset X$ be a domain, and $C_{bs}\geq 1$ be a constant. We say that $D$ satisfies the {\it $C_{bs}$-ball separation condition}, if for all $x$, $y$ in $D$ and for each
quasihyperbolic geodesic $\gamma$ joining $x$ and $y$, we have for every $z\in \gamma$,
$$B(z,C_{bs}d_D(z)) \cap \beta_{x,y} \not=\emptyset ,$$
where $\beta_{x,y}$ is any other curve joining $x$ and $y$ in $D$.
\edefe

\bdefe
Let $(X,d)$ be a metric space, $D\subset X$ be a domain and let $c\geq 1$ be a constant. We say that $D$ is
\begin{enumerate}  \item {\it $c$-$LLC_1$}, if for all $x\in D$ and $r>0$, we have every pair of points in $B(x,r)$ can be joined by a curve in $B(x,cr)$.

  \item {\it $c$-$LLC_2$}, if for all $x\in D$ and $r>0$, we have every pair of points in $D\backslash B(x,r)$ can be joined by a curve in  $D\backslash B(x,\frac{r}{c})$.
 \item {\it $c$-$LLC$}, if it is both $c$-$LLC_1$ and $c$-$LLC_2$.
\end{enumerate}
Moreover, $D$ is called {\it linearly locally connected} or {$LLC$}, if it is $c$-$LLC$ for some constant $c\geq 1$.
\edefe

\bdefe Let $c\geq 1$. A noncomplete metric space $(X,d)$ is {\it $c$-locally externally connected} ($c$-$LEC$) provided the $c$-$LLC_2$ property holds for all points $x\in X$ and all $r\in (0,d(x)/c)$.
\edefe

In \cite{BH}, Buckley and Herron obtained the following interesting characterization of uniform metric spaces.

\begin{Thm}\label{bhthm4.2}$($\cite[Theorem 4.2]{BH}$)$ A minimally nice metric space $(X, d)$ is uniform and $LEC$ if and
only if it is quasiconvex, $LLC$ with respect to curves, and satisfies a weak slice condition.
These implications are quantitative.\end{Thm}

\bdefe A metric space $(X,d)$ is called {\it annular quasiconvex}, if there is a constant $\lambda \geq
1$ so that, for any $x\in X$ and all $0 < r' < r$, each pair of points $y, z$ in
$B(x, r) \backslash B(x, r')$ can be joined with a curve $\gamma_{yz}$ in $B(x, \lambda r)\backslash B(x, r'/\lambda)$
such that $\ell(\gamma_{yz})\leq \lambda d(y, z)$.
\edefe

It is not difficult to see that $\lambda$-annularly quasiconvexity property implies $C$-$LLC_2$, and hence $C$-$LEC$, where $C=2\lambda^2$.
\medskip

\subsection{The proof of Theorem \ref{thm-2}}

Necessity: Suppose that $D$ is uniform. Then we know that $D$ is John and quasiconvex. Moreover, it follows from \cite[Theorem 3.6]{BHK} that $(D,k)$ is a roughly starlike Gromov hyperbolic space because $D$ is bounded, where $k$ is the quasihyperbolic metric of $D$. It remains to show that $D$ is $LLC$. Since $X$ is $A$-annularly quasiconvex, it follows that $D$ is $LEC$. Then we deduce from Theorem \Ref{bhthm4.2} that $D$ is $LLC$.

Sufficiency: To prove the uniformity of $D$, we only need to prove that every quasihyperbolic geodesic $\gamma$ in $D$ is a uniform arc. We assume that $D$ is $c$-quasiconvex and $\delta$-hyperbolic. By \cite[Theorem 1.2]{KLM14}, we find that $D$ satisfies both the $C_{gh}$-Gehring-Hayman condition and the $C_{bs}$-ball condition for some constants $C_{gh},C_{bs}\geq 1$. So to prove the sufficiency, we only need to show that each quasihyperbolic geodesic in $D$ is a cone arc.

We first assume that $D$ is $a$-John. Since $D$ is a bounded $\delta$-hyperbolic domain of $X$, we see from \cite[Theorem 3.1]{BB03} that $(D,k)$ is $K$-roughly starlike, because $X$ is annularly quasiconvex. Then from Theorem \ref{thm-1} the uniformity of $D$ follows.

We are thus left to assume that $D$ is $c_0$-LLC. Again by virtue of the Gehring-Hayman condition, we only need to show that there is a uniform upper bound for the constant $\Lambda$ such that
$$\min\{\ell(\gamma[x,z]),\ell(\gamma[z,y])\}=\Lambda d_D(z)$$
for each pair of points $x,y\in D$, for any quasihyperbolic geodesic $\gamma$ in $D$ joining $x$ and $y$, and for every point $z\in \gamma$.

To this end, we deduce from the $C_{gh}$-Gehring-Hayman condition that $$\ell(\gamma[x,z])\leq cC_{gh}d(x,z)\;\;{\rm and}\;\;\ell(\gamma[y,z])\leq cC_{gh}d(y,z),$$
because the subarcs $\gamma[x,z]$ and $\gamma[x,y]$ are also quasihyperbolic geodesics.

Thus we have
$$\min\{d(x,z),d(y,z)\}\geq \frac{\Lambda}{cC_{gh}}d_D(z).$$
On the other hand, since $D$ is $c_0$-LLC, we know that there is a curve $\beta$ joining $x$ to $y$ with
\beq\label{eq-new1}\beta\subset X\setminus \overline{B}(z,\frac{\Lambda}{cc_0C_{gh}}d_D(z)).\eeq
Furthermore, since $\gamma$ is a quasihyperbolic geodesic and $D$ satisfies the $C_{bs}$-ball separation condition, we see that $$\beta\cap B(z,C_{bs}d_D(z))\not=\emptyset,$$ which together with \eqref{eq-new1} show that
$$\Lambda\leq cc_0C_{gh}C_{bs},$$
as required.

Hence, the proof of Theorem \ref{thm-2} is complete.


\end{document}